\documentclass[11pt, oneside]{article}
\usepackage{amsmath,amsfonts,amssymb,amscd}
\usepackage{amsthm}
\usepackage{eucal}
\usepackage{enumerate}
\date{}

\begin{document}

\bibliographystyle{plain}

\newtheorem{theorem}{Theorem}[section]
\newtheorem{proposition}{Proposition}[section]
\newtheorem{lemma}{Lemma}[section]
\newtheorem{corollary}[lemma]{Corollary}

\title{Manifolds of positive Ricci curvature, quadratically asymptotically nonnegative curvature, and infinite Betti numbers}

\author{Huihong Jiang and Yi-Hu Yang\thanks{Partially supported by NSF of China (No.11571228)}}

\maketitle

\begin{abstract}
In a previous paper \cite{JY}, we constructed complete manifolds of positive Ricci curvature
with quadratically asymptotically nonnegative curvature and infinite topological type
but dimension $\ge 6$. The purpose of the present paper is to use a different technique to exhibit
a family of complete $I$-dimensinal ($I\ge5$) Riemannian manifolds of positive Ricci curvature,
quadratically asymptotically nonnegative sectional curvature,
and certain infinite Betti number $b_j$ ($2\le j\le I-2$).
\end{abstract}

\section{Introduction}
An interesting topic in Riemannian geometry is to give some topological constraints of complete
manifolds under certain prescribed curvature assumptions. For example, a remarkable theorem
due to Gromov says that the total Betti number of complete manifolds $M^n$ (either compact or noncompact)
with nonnegative sectional curvature is bounded by a constant only depending on $n$ \cite{G}.
A natural question is: {\it Can one bound the Betti numbers of open manifolds with nonnegative Ricci curvature?}

For the first Betti number, Anderson \cite{And} proved that $b_1(M^n)\le n$ for a complete manifold
with nonnegative Ricci curvature and $b_1(M^n)\le n-3$ if the manifold has positive Ricci curvature.
For the codimension one Betti number of open manifolds with nonnegative Ricci curvature, Shen and Sormani \cite{SSo}
showed that either $M^n$ is a flat normal bundle over a compact totally geodesic submanifold or $M^n$
has a trivial codimension one integer homology. So, $b_{n-1}(M^n)={\text{dim}}H_{n-1}(M^n;\mathbb{Z})\le1$.

For other Betti numbers $b_k(M^n)$, $2\le k\le n-2$, it is difficult to find out some similar control,
even with some additional assumptions. In fact, Sha and Yang \cite{ShY1} (also cf. \cite{ShY,SW})
constructed $n$-dim ($n\ge4$) open manifolds with positive Ricci curvature and certain infinite Betti number
$b_k$, $2\le k\le n-2$. Later, using a totally different method, Menguy \cite{M1,M2} also constructed
$4$-dim open manifolds with positive Ricci curvature, maximal volume growth or minimal volume growth,
and infinite second Betti number; moreover, he showed that similar things can be done so that one can get
higher dimensional examples with infinite even dimensional Betti numbers $b_k$, $2\le k\le n-2$.

So, in order to get certain control of Betti numbers, it seems to needs some additional assumptions.
There is a natural and interesting sectional curvature condition, {\it quadratically asymptotically
nonnegative curvature}, which can be defined as follows. Let
$$
K_{p_0}(t)=\inf_{M^n\setminus B(p_0,t)}K,
$$
where $K$ denotes the sectional curvature of $M$, and the infimum is taken over all the 2-planes at the points
in $M\setminus B(p_0,t)$. We say that $M$ is of quadratically asymptotically nonnegative curvature
if for all $t\geq 0$,
$$
K_{p_0}(t)\geq-\frac{K_0}{1+t^2},
$$
for some positive constant $K_0$. We also say that $M$ is of finite topological type if it is homeomorphic to
the interior of a compact manifold with boundary, otherwise, infinite topological type. In particular,
manifolds of finite topological type must have finite Betti numbers. A little additional computation
shows that all the examples mentioned above are not of quadratically asymptotically nonnegative curvature.

Actually, Sha and Shen \cite{SS} conjectured that a complete manifold $M^n$ of nonnegative Ricci curvature
and quadratically asymptotically nonnegative curvature should be of finite topological type.
In \cite{JY}, the authors constructed a counterexample with infinite second Betti number
in the case of dimension not less than $6$, based on works of Perelman \cite{P} and Menguy \cite{M1, M2, M3}.

The purpose of the present paper is to use a topologically $\&$ metrically different technique to construct
a family of complete $I$-dim ($I\ge5$) Riemannian manifolds of positive Ricci curvature,
quadratically asymptotically nonnegative sectional curvature, and certain infinite Betti number $b_j$
($2\le j\le I-2$). In particular, this again gives a negative answer of Sha and Shen's conjecture in dimension $5$.
More precisely, we have the following

\begin{theorem} For all integers $m,n$ satisfying $n\ge m\ge2$ and $n\ge3$, there exists an $(m+n)$-dimensional
complete Riemnannian manifold $M$, which is of positive Ricci curvature, quadratically asymptotically
nonnegative curvature, and infinite Betti numbers $b_m, b_n$. In particular, it is of infinite topological type.
\end{theorem}

\noindent
{\bf Remarks:}
\begin{enumerate}[1)]
  \item U. Abresch showed in \cite{A} that if a complete $n$-dim ($n\ge3$) manifold satisfies
  that the integral $\int_0^{\infty}rK^-_{p_0}(r)dr=b_0<\infty$, here $K^-_{p_{0}}(r)=\max\{-K_{p_{0}} (r), 0\}$
  (a certain kind of decay on the lower sectional curvature bound which is a little stronger
  than quadratically asymptotically nonnegative curvature), then the manifold is of finite
  toplogical type; moreover, the total Betti number of the manifold is uniformly bounded.

  \item For the 3-dim case, complete noncompact Riemannian manifolds with nonnegative Ricci curvature
  were completely classified by G. Liu \cite{L} based on the technique earlier developed
  by Schoen-Yau \cite{SY}, which must be of finite topological type.

  \item Combined our examples in dim $\ge5$ with Liu's result in dim $3$, the remained is the case of dim $4$.
  But, techniquely, it is difficult to use the present method to construct counterexamples of dim $4$.
  On the other hand, this seems to be the most interesting case in some sense.
\end{enumerate}
\vskip .2cm

Topologically, the construction of the manifolds is similar to that of Sha and Yang \cite{ShY1},
i.e. removing some balls $D_i^{m+1}$, $1\le i<\infty$, from $\mathbb{R}^{m+1}$,
and then gluing $(\mathbb{R}^{m+1}\setminus D_i^{m+1})\times S^{n-1}$ with $P_i=S^m\times D^n$ along the boundary,
which is different from our previous one \cite{JY}. Then, the resulting $(m+n)$-dim manifold has infinite Betti
numbers $b_m$ and $b_n$.

The construction for the metrics however is totally different from that of Sha and Yang; in their construction,
a $C^1$ metric on the whole manifold is directly given. Instead, we'll first give a $C^0$ metric on the whole manifold;
more precisely, we will give a $C^1$ metric on each part 
(i.e. $(\mathbb{R}^{m+1}\setminus D_i^{m+1})\times S^{n-1}$ and $P_i=S^m\times D^n$) but with isometric boundaries and then glue them
together along the boundaries, so that the metric is only $C^0$ on the glued part, but $C^1$ on the remained part.
On the other hand, the above $C^1$ metrics actually have some additional restrictions on the normal curvatures
of the corresponding boundaries, so that we next can use a $C^1$ metric to replace
the $C^0$ metric near the glued part, moreover the corresponding Ricci and sectional curvature
properties are preserved. Thus, we can get a required complete $C^1$ metric on the whole manifold.
Essentially, the existence to such a metric is guaranteed by the following gluing criterion due to G. Perelman
(cf. \cite{P, BWW}, also see the remark in \S 4);
but here for the sake of clarity, we will give an explicit $C^1$ construction.

\vskip .2cm
\noindent {\bf Gluing Criterion:} {\it Let $M_{1}$, $M_{2}$ be two compact smooth manifolds
of positive Ricci curvature, with isometric boundaries $\partial M_{1} \backsimeq \partial M_{2} = X$.
Suppose that the normal curvatures of $\partial M_{1}$ are bigger than the negative of
the normal curvatures of $\partial M_{2}$. Then, the glued manifold $M_{1}\cup_{X} M_{2}$
along the boundary $X$ can be smoothed near $X$ to produce a manifold of positive Ricci curvature.}
\vskip .2cm

Based the above idea, the construction can be divided into the following three steps. First, we will construct a $C^1$ metric $ds^2$
on $Q=\mathbb{R}^{m+1}\times S^{n-1}=[t_0,+\infty)\times_{u(t)}S^m\times_{g(t)}S^{n-1}$ and
a $C^1$ metric $d\bar{s}^2$ on $P_i=S^m\times D^n$ ($1\le i<\infty$). These are done in \S 2.

Then, we will compute the curvature tensors and show that both metrics have positive Ricci curvature
and quadratically asymptotically nonnegative sectional curvature by choosing some appropriate constants.
These will be done in \S 3.

Finally, we will remove the geodesic balls $B_{\frac45r_i}(o_i)$ from $[t_0,+\infty)\times_{u(t)}S^m$
(for the precise meaning of $r_i$ and $o_i$, see \S 2.1) and glue
$Q\setminus\coprod\limits_{i=1}^{+\infty}(B_{\frac45r_i}(o_i) \times_{g_i}S^{n-1})$
with $P_i$ along the boundary $\partial B_{\frac45r_i}(o_i)\times_{g_i}S^{n-1}=S^m\times S^{n-1}$;
this will give a $C^0$ metric $h_0$ on the whole manifold. Then, we will use a $C^1$ metric $h_1$ to replace the above $h_0$
near the glued boundary and verify that it still has positive Ricci curvature and quadratically
asymptotically nonnegative curvature. These will be done in \S 4.

So, our manifold is actually
\begin{align*}
M^{m+n}&=\left(Q\setminus\coprod\limits_{i=1}^{+\infty}(B_{\frac45r_i}(o_i)\times_{g_i}S^{n-1})\right)
\cup_{\textrm{Id}}\coprod\limits_{i=1}^{+\infty}P_i\\
&=\left((\mathbb{R}^{m+1}\setminus\coprod\limits_{i=1}^{+\infty}D_i^{m+1})\times S^{n-1}\right)
\cup_{\textrm{Id}}\coprod\limits_{i=1}^{+\infty}(S^m\times D^n)_i,
\end{align*}
where by 'Id' we mean gluing along the corresponding boundaries through the identity map.

\section{Construction of the metrics}
As seen in the introduction, our topological manifold is the following
\begin{align*}
M^{m+n}&=\left(Q\setminus\coprod\limits_{i=1}^{+\infty}(B_{\frac45r_i}(o_i)\times_{g_i}S^{n-1})\right)
\cup_{\textrm{Id}}\coprod\limits_{i=1}^{+\infty}P_i.
\end{align*}
In the following, we'll construct the required metrics on $Q$ and $P_i$ ($i=1, 2, \cdots$) respectively.

\subsection{Construction of $Q$}
We equip $Q=[t_0,+\infty)\times_{u(t)}S^m\times_{g(t)}S^{n-1}$ (for some $t_0=t_1-\psi_1$, $t_1$
and $\psi_1$ will be given later) with the following metric
$$
ds^2=dt^2+u^2(t)d\theta_{S^m}^2+g^2(t)d\sigma_{S^{n-1}}^2,
$$
where $d\theta_{S^m}^2$ and $d\sigma_{S^{n-1}}^2$ are the standard metrics on $S^m$ and $S^{n-1}$ respectively.

Give some constants $c, \gamma, \alpha$ and $t_1$ satisfying
$$
0<c<\frac13,~~0<\gamma<\frac14,~~\alpha>1,~~t_1> 1.
$$
Set
$$
K=K(c)=\frac{1-c^2}{c^2},
$$
and define $\psi=\psi(c)$ by
$$
\sin(\sqrt{K}\psi)=\sqrt{1-c^2}.
$$
Take $r>0$ satisfying
$$
r \le r(c)=\frac{\pi}{4\sqrt{K}}-\frac{\psi}2.
$$
Note that $c$ first is fixed while $r$, $\gamma, \alpha$ and $t_{1}$ will be determined later.
We define, for $i=1, 2, \cdots,$
$$
t_i=t_1\alpha^{i-1},
$$
$$
r_i=rt_i,
$$
$$
\psi_i=\psi t_i,
$$
$$
K_i=\frac K{t_i^2},
$$
and
$$
\Delta=\sqrt{K_i}(2r_i+\psi_i)=\sqrt{K}(2r+\psi).
$$

In the following, we'll give the constructions of $u$ and $g$ respectively.

\subsubsection{Construction of $u(t)$}
For $t_0\le t\le t_1$, we set
$$
u(t)=\frac1{\sqrt{K_1}}\sin\left(\sqrt{K_1}(t-t_1+\psi_1)\right),
$$
Here, $t_0=t_1-\psi_1=(1-\psi)t_1$.

For $t_i<t<t_i+2r_i$, $i=1, 2, \cdots$, we set
$$
u(t)=\frac1{\sqrt{K_i}}\sin\left(\sqrt{K_i}(t-t_i+\psi_i)\right);
$$
for $t_i+2r_i<t<t_{i+1}$, we set
$$
u(t)=t_i w(\frac t{t_i}),
$$
where $w$ is a $C^2$ function on $[1+2r,\alpha]$, which is independent of $i$.

Actually, we can define $w$ as follows. First, set
$$
w(t)=\min\{\frac{\sin\Delta}{\sqrt{K}}+(t-1-2r)\cos\Delta+
\frac{c+1}{2\log\frac{\alpha}{1+2r}}(t\log\frac{t}{1+2r}-t+1+2r), ct \},
$$
and then smoothen $w(t)$ to be a $C^2$ function. Since, at $t=1+2r$,
$$
\frac{\sin\Delta}{\sqrt{K}}<c(1+2r);
$$
while at $t=\alpha$,
$$
\frac{\sin\Delta}{\sqrt{K}}+(\alpha-1-2r)\cos\Delta+
\frac{c+1}{2\log\frac{\alpha}{1+2r}}(\alpha\log\frac{\alpha}{1+2r}-\alpha+1+2r)>c\alpha,
$$
provided with $\cos\Delta\ge\frac{c+1}{2\log\frac{\alpha}{1+2r}}$, which is possible for $\alpha\ge\alpha_0(c,r)$.

Thus, $u(t)$ is $C^1$ at the endpoints $t_i+2r_i=(1+2r)t_i$ and $t_{i+1}=\alpha t_i$. On the other hand, if $w(t)=ct$,
$$
\begin{cases}
w_t \equiv c,\vspace{+0.2cm}\\
w_{tt} \equiv 0;
\end{cases}
$$
and if $w(t)=\frac{\sin\Delta}{\sqrt{K}}+(t-1-2r)\cos\Delta+
\frac{c+1}{2\log\frac{\alpha}{1+2r}}(t\log\frac{t}{1+2r}-t+1+2r)$,
$$
\begin{cases}
\cos\Delta \le w_t \le \frac{1+3c}2 < 1,\vspace{+0.2cm}\\
-3\frac{w_{tt}}w+\frac{\gamma(1-2\gamma)}{t^2} > 0,
\end{cases}
$$
when $\cos\Delta\ge\frac{c+1}{\log\frac{\alpha}{1+2r}}$ and $\alpha\ge\alpha_{1}(c,r,\gamma)$.

\vskip .2cm

\noindent
{\bf Conclusion:} When $\alpha\ge\alpha_{2}(c,r,\gamma)$, we have $u(t)$ satisfying, for $t>t_{1}$,
$$
\begin{cases}
\cos\Delta \le u_t \le c, \quad -\frac{u_{tt}}u=\frac{K}{t_i^2}, & t_i<t<t_i+2r_i,\vspace{+0.2cm}\\
\cos\Delta \le u_t \le \frac{1+3c}2, \quad -3\frac{u_{tt}}u+\frac{\gamma(1-2\gamma)}{t^2} > 0, & t_i+2r_i<t< t_{i+1}.
\end{cases}
$$

\subsubsection{Construction of $g(t)$}
For $t_0\le t\le t_1+\frac{r_1}6$, we set
$$
g(t)\equiv g_1=g(t_1+\frac{r_1}6)=(t_1+\frac{r_1}6)^{\gamma}=(1+\frac r6)^{\gamma}t_1^{\gamma}.
$$
For $t_i+\frac{r_i}6<t<t_{i+1}+\frac{r_{i+1}}6=\alpha(t_i+\frac{r_i}6)$, $i=1, 2, \cdots$, we set
$$
\frac{g_t}g=\begin{cases}
0, & t_i+\frac{r_i}6<t<t_i+\frac{11r_i}6,\vspace{+0.2cm}\\
\frac{6\gamma\beta}{(1+2r)rt_i^2}\left(t-(t_i+\frac{11r_i}6)\right),&t_i+\frac{11r_i}6<t<t_i+2r_i,\vspace{+0.2cm}\\
\frac{\gamma\beta}t,&t_i+2r_i<t<t_{i+1},\vspace{+0.2cm}\\
\frac{6\gamma\beta}{\alpha^2rt_i^2}\left((t_{i+1}+\frac{r_{i+1}}6)-t\right),&t_{i+1}<t<t_{i+1}+\frac{r_{i+1}}6,
\end{cases}
$$
where $\beta=\frac{\log\alpha}{\log\alpha-\log(1+2r)+\frac{r(1+r)}{6(1+2r)}}=\beta(r,\alpha)$, so that
$$
g(t_i+\frac{r_i}6)=(t_i+\frac{r_i}6)^{\gamma}, \quad \forall i\ge1.
$$
Thus, for $t_i+\frac{r_i}6<t<t_{i+1}+\frac{r_{i+1}}6=\alpha(t_i+\frac{r_i}6)$,
$$
g(t)\le g(t_{i+1}+\frac{r_{i+1}}6)=\alpha^{\gamma}(t_i+\frac{r_i}6)^{\gamma}\le\alpha^{\gamma}t^{\gamma},
$$
and
$$
g(t)\ge g(t_i+\frac{r_i}6)=(t_i+\frac{r_i}6)^{\gamma}\ge\alpha^{-\gamma}t^{\gamma}.
$$
Note that for $\alpha\ge\alpha_3(r)$, we can have $\frac12\le\beta\le2$. Then, one has
$$
|\frac{g_t}g|\le\frac{2(1+\frac r6)\gamma}t
$$
and
$$
\begin{cases}
|\frac{g_{tt}}g|\le\frac{12\gamma(1+\frac{r\gamma}3)}{rt_i^2},&t_i<t<t_i+2r_i,\vspace{+0.2cm}\\
-\frac{g_{tt}}g\ge\frac{\gamma(1-2\gamma)}{2t^2},&t_i+2r_i<t<t_{i+1}.
\end{cases}
$$

\vskip .2cm

\noindent
{\bf Conclusion:} When $\alpha\ge\alpha_3(r)$, we have a $C^1$ function $g(t)$ satisfying,  for $t>t_1$,
$$
\alpha^{-\gamma}t^{\gamma}\le g(t)\le\alpha^{\gamma}t^{\gamma},
$$
$$
|\frac{g_t}g|\le\frac{2(1+\frac r6)\gamma}t,
$$
and
$$
\begin{cases}
|\frac{g_{tt}}g|\le\frac{12\gamma(1+\frac{r\gamma}{3})}{rt_i^2},&t_i<t<t_i+2r_i,\vspace{+0.2cm}\\
-\frac{g_{tt}}g\ge\frac{\gamma(1-2\gamma)}{2t^2},&t_i+2r_i<t<t_{i+1}.
\end{cases}
$$

\vskip .2cm
\noindent {\bf Remark:} We here remark that, at some discrete points, $u(t)$ and $g(t)$ are only $C^1$.
However, if the manifold constructed in this manner has positive Ricci curvature on the complement
of those $C^1$ parts, the manifold can then be smoothen to be a $C^2$ manifold of positive Ricci curvature;
for this, one can refer to \cite{P}, also \cite{BWW, M1}.

\subsection{Construction of $P_i$ ($1\leq i<+\infty$)}
$P_i$ topologically is $S^m \times D^n$ but with a metric
$$
d\bar{s}^2=d\bar{t}^2+\bar{u}^2(\bar{t})d\theta_{S^m}^2+\bar{g}^2(\bar{t})d\sigma_{S^{n-1}}^2
$$
where $d\theta_{S^m}^2$ and $d\sigma_{S^{n-1}}^2$ are the standard metrics on $S^m$ and $S^{n-1}$ respectively,
and $0\le\bar{t}\le R_i$ ($R_i$ will be fixed in the following).

First, some constants $C_1, C_2, C_3, N, R_i, b_i, l_i, C_4$ are given as follows,
$$
C_1=(\frac{\sin\frac45\sqrt{K}r}{\sqrt{K}})^{-1},
$$
$$
C_2=\cos\frac35\sqrt{K}r>\cos\frac45\sqrt{K}r,
$$
$$
0<C_3<1-C_2^2,
$$
$$
N\ge5n,
$$
$$
R_i>\frac{C_4}{C_1C_2(1-C_3)}t_i,
$$
$$
b_i=R_i-\frac{C_4}{C_1C_2(1-C_3)}t_i,
$$
$$
l_i=\frac{C_1(NC_2)^{1/C_3}}{t_i},
$$
and $C_4$ is given by $C_2(1-C_4)^{\frac{C_3}{1-C_3}}={\frac 1N}$, provided with $N>\frac{1}{C_2}$;
clearly, $0<C_4<1$. Note that $C_1,C_2$ are given by $c, r$; $C_3$, $N$ and $R_i$ will be fixed later;
while $b_i, l_i$ and $C_4$ are given by $C_1,C_2,C_3,N$, and $R_i$.

\subsubsection{Construction of $\bar{u}(\bar{t})$}
$\bar{u}(\bar{t})$ is given by
$$
\frac{\bar{u}_{\bar{t}}}{\bar{u}}=\begin{cases}
\frac{1}{(1-C_3)\bar{t}+\frac{N}{l_i}-(1-C_3)b_i}, &b_i\le\bar{t}\le R_i,\vspace{+0.2cm}\\
\frac{l_i}{Nb_i}\bar{t}, &0\le\bar{t}\le b_i,
\end{cases}
$$
and $\bar{u}(R_i)=\frac{t_i}{C_1}, \bar{u}_{\bar{t}}(R_i)=C_2$. Thus,
$$
\frac{\bar{u}_{\bar{t}\bar{t}}}{\bar{u}}=\begin{cases}
C_3\frac{\bar{u}^2_{\bar{t}}}{\bar{u}^2}
=\frac{C_3}{\left((1-C_3)\bar{t}+\frac{N}{l_i}-(1-C_3)b_i\right)^2}, &b_i\le\bar{t}\le R_i,\vspace{+0.2cm}\\
\frac{l_i}{Nb_i}+(\frac{l_i}{Nb_i}\bar{t})^2, &0\le\bar{t}\le b_i.
\end{cases}
$$
Then, for $b_i\le\bar{t}\le R_i$,
$$
\bar{u}(\bar{t})
=C_2(\frac{t_i}{C_1C_2})^{-\frac{C_3}{1-C_3}}\big((1-C_3)\bar{t}+\frac{N}{l_i}-(1-C_3)b_i\big)^{\frac{1}{1-C_3}},
$$
and $\bar{u}(b_i)=\frac{1}{l_i}, \bar{u}_{\bar{t}}(b_i)=\frac1N$; while for $0\le\bar{t}\le b_i$,
$$
0<\frac{\bar{u}_{\bar{t}\bar{t}}}{\bar{u}}\le\frac{l_i}{Nb_i}+(\frac{l_i^2}{N^2})^2,
$$
and $\bar{u}(0)=\frac{1}{l_i}\exp(-\frac{l_ib_i}{2N})>0$.

\subsubsection{construction of $\bar{g}(\bar{t})$}
Setting
$$
A=\frac{(1+\frac{r}{6})^\gamma}{\left(\frac{1}{C_1(NC_2)^{1/{C_3}}}+\frac{C_4}{C_1C_2(1-C_3)}\right)^\gamma}
=A(r,c,\gamma,N,C_3)=O(1),
$$
$d_i=d_i(A,t_i)$ then is given by the equation
$$
d_i^{1-\gamma}\cos(l_id_i)=\gamma A.
$$
We remark that when $d_i=0$, the left hand side of the equation is 0, which is smaller than
the right hand side $\gamma A$; and when $d_i=O(t_i)$, the left hand side is $O(t_i^{1-\gamma})$
(since $l_i=O(\frac{1}{t_i})$), which is much bigger than the right hand side.
Thus, there exists such a positive $d_i$ satisfying the above equation; moveover, $d_i=O(1)$.

Set $B=\frac{\sin(l_id_i)}{l_i}-Ad_i^\gamma=Ad_i^\gamma(\gamma\frac{\tan(l_id_i)}{l_id_i}-1)$,
which is negative when $t_i$ is sufficiently large (i.e. $i$ is sufficiently large). Moreover, $B=O(1)$.

Finally, $R_i$ is given by the following equation
$$
AR_i^\gamma+B=(1+\frac{r}{6})^\gamma t_i^\gamma,
$$
i.e.
$$
\left(b_i+\frac{C_4}{C_1C_2(1-C_3)}t_i\right)^\gamma=R_i^\gamma=\frac{(1+\frac{r}{6})^\gamma t_i^\gamma-B}{A}.
$$
Since $B<0$ and $B=O(1)$, we have
\begin{align*}
\frac{(1+\frac{r}{6})^\gamma t_i^\gamma}{A}
&\le\frac{(1+\frac{r}{6})^\gamma t_i^\gamma-B}{A}\\
&\le\frac{(1+\frac{r}{6})^\gamma t_i^\gamma}{A}\left(\frac{\frac{2}{C_1(NC_2)^{1/{C_3}}}+\frac{C_4}{C_1C_2(1-C_3)}}
{\frac{1}{C_1(NC_2)^{1/{C_3}}}+\frac{C_4}{C_1C_2(1-C_3)}}\right)^\gamma,
\end{align*}
as $t_i>t'(c,r,\gamma,C_3,N)$, i.e. $t_i$ sufficiently large.

Then, we get
$$
R_i>\frac{C_4}{C_1C_2(1-C_3)}t_i
$$
and
$$
\frac 1{l_i}\le b_i\le\frac 2{l_i}.
$$
Moreover, $R_i=O(t_i)$ and $b_i=O(t_i)>d_i$, since $l_i=O(\frac 1{t_i})$.

Now, $\bar{g}(\bar{t})$ can be defined as follows
$$
\bar{g}(\bar{t})=\begin{cases}
A\bar{t}^\gamma+B, &d_i\le\bar{t}\le R_i,\vspace{+0.2cm}\\
\frac{\sin(l_i\bar{t})}{l_i}, &0\le\bar{t}\le d_i,
\end{cases}
$$
which satisfies
$$
\frac{\bar{g}_{\bar{t}}}{\bar{g}}=\begin{cases}
\frac{\gamma A\bar{t}^{\gamma-1}}{A\bar{t}^\gamma+B}, &d_i\le\bar{t}\le R_i,\vspace{+0.2cm}\\
\frac{l_i}{\tan(l_i\bar{t})}, &0\le\bar{t}\le d_i
\end{cases}
$$
and
$$
\frac{\bar{g}_{\bar{t}\bar{t}}}{\bar{g}}=\begin{cases}
-\frac{\gamma(1-\gamma)A\bar{t}^{\gamma-2}}{A\bar{t}^\gamma+B}, &d_i\le\bar{t}\le R_i,\vspace{+0.2cm}\\
-l_i^2, &0\le\bar{t}\le d_i.
\end{cases}
$$

\section{Quadratically asymptotic non-negativeness of curvature and positiveness of Ricci curvature}
In this section, we will calculate the curvature tensors of $Q$ and $P_i$ with the given metrics in \S 2,
and show that they have both positive Ricci curvature and quadratically asymptotically nonnegative
sectional curvature by taking some appropriate constants, i.e. $\eta, r, \gamma, C_3$
sufficiently small and $\alpha, t_1, N$ sufficiently big.

\subsection{Curvatures of $Q$}
For $Q$, the metric is
$$
ds^2=dt^2+u^2(t)d\theta_{S^m}^2+g^2(t)d\sigma_{S^{n-1}}^2,
$$
where $m\ge2, n\ge3$ and $n\ge m$, thus $2(n-1)\ge m\ge2$. Let $T, \{\Theta_k\}, \{\Sigma_l\}$ be an
orthonormal basis of the tangent space corresponding to the directions
$dt, d\theta^2_{S^m}, d\sigma^2_{S^{n-1}}$ respectively.
Then, for $t\ge t_1$, the sectional curvatures can be computed as follows.

\begin{align*}
K(T,\Theta_k,\Theta_k,T)&=-\frac{u_{tt}}{u}
\ge-\frac{\gamma(1-2\gamma)}{3t^2},\\
K(T,\Sigma_l,\Sigma_l,T)&=-\frac{g_{tt}}{g}
\ge-\frac{12\gamma(1+\frac{r\gamma}{3})(1+2r)}{rt^2},\\
K(\Theta_k,\Theta_p,\Theta_p,\Theta_k)&=\frac{1}{u^2}-\frac{u_t^2}{u^2}
\ge\frac{1-(\frac{1+3c}{2})^2}{\cos^2\Delta}\frac{1}{t^2}>0,\\
K(\Sigma_l,\Sigma_q,\Sigma_q,\Sigma_l)&=\frac{1}{g^2}-\frac{g_t^2}{g^2}
\ge\frac{1}{t^{2\gamma}}-\frac{4(1+\frac r6)^2\gamma^2}{t^2}>0,\\
K(\Theta_k,\Sigma_l,\Sigma_l,\Theta_k)&=-\frac{u_t}{u}\frac{g_t}{g}
\ge-\frac{(1+3c)(1+\frac r6)\gamma}{\cos\Delta}\frac{1}{t^{2}};
\end{align*}
and other terms of curvature tensors are zero. Thus, we have
$$
K_{p_0}(t)\ge-\frac{K_0^2(c,r,\gamma)}{t^2}.
$$
i.e. $Q$ has quadratically asymptotically nonnegative sectional curvatures.

The nonzero Ricci curvatures of $Q$ are as follows.
\begin{align*}
Ric(T,T)&=-m\frac{u_{tt}}{u}-(n-1)\frac{g_{tt}}{g}\\
        &\ge\begin{cases}
         -m\frac{u_{tt}}{u}+(n-1)\frac{\gamma(1-2\gamma)}{2t^2}, &t_i+2r_i<t<t_{i+1},\vspace{+0.2cm}\\
         m\frac{K}{t_i^2}-(n-1)\frac{12\gamma(1+\frac{r\gamma}{3})}{rt_i^2}, &t_i<t<t_i+2r_i,
         \end{cases}\\
         &\ge\begin{cases}
         \frac m3(-3\frac{u_{tt}}{u}+\frac{\gamma(1-2\gamma)}{t^2})>0, &t_i+2r_i<t<t_{i+1},\vspace{+0.2cm}\\
         m\frac{K}{t_i^2}-(n-1)\frac{12\gamma(1+\frac{r\gamma}{3})}{rt_i^2}>0, &t_i<t<t_i+2r_i,
         \end{cases}
\end{align*}
for $i=1, 2, \cdots$, and provided with $\gamma<\gamma_0(n,c,r)$;

\begin{align*}
Ric(\Theta_k,\Theta_k)&=(m-1)(\frac1{u^2}-\frac{u_t^2}{u^2})-\frac{u_{tt}}u-(n-1)\frac{u_t}u\frac{g_t}g\\
         &\ge(m-1)(\frac1{u^2}-\frac{u_t^2}{u^2})-\frac13\frac{\gamma(1-2\gamma)}{t^2}
         -(n-1)\frac{u_t}{u}\frac{2(1+\frac r6)\gamma}t\\
         &=\frac1{u^2}\left((m-1)-(m-1)u_t^2-2(n-1)\frac{(1+\frac r6)\gamma u_tu}t\right)
         -\frac13\frac{\gamma(1-2\gamma)}{t^2}\\
         &\ge\frac1{(\frac{1+3c}2)^2t^2}\left((m-1)(1-(\frac{1+3c}2)^2)-2(n-1)(1+ \frac r6)\gamma(\frac{1+3c}2)^2\right)\\
         &\quad-\frac13\frac{\gamma(1-2\gamma)}{t^2}\\
         &=\frac1{(\frac{1+3c}2)^2t^2}\Big((m-1)-(\frac{1+3c}2)^2((m-1)+2(n-1)(1+\frac r6)\gamma\\
         &\quad+\frac{\gamma(1-2\gamma)}3)\Big)\\
         &>0,
\end{align*}
provided with $\gamma<\gamma_1(m,n,c,r)$; while

\begin{align*}
Ric(\Sigma_l,\Sigma_l)&=(n-2)(\frac1{g^2}-\frac{g_t^2}{g^2})-\frac{g_{tt}}g-m\frac{u_t}u\frac{g_t}g\\
         &\ge(n-2)\left(\frac1{\alpha^{\gamma}t^{\gamma}}-\frac{4(1+\frac r6)^2\gamma^2}{t^2}\right)
         -\frac{12\gamma(1+2r)(1+\frac{r\gamma}3)}{rt^2}\\
         &\quad-m\frac{(1+3c)(1+\frac r6)\gamma}{t^2\cos\Delta}\\
         &>0,
\end{align*}
provided with $t_1>t_1(c,r,\gamma)$. Moreover, all the off-diagonal terms of the Ricci tensor vanish.

On the other hand, when $t_0\le t\le t_1$, the metric is actually
$$
ds^2=dt^2+\left(\frac1{\sqrt{K_1}}\sin\left(\sqrt{K_1}(t-t_1+\psi_1)\right)\right)^2d\theta_{S^m}^2+g_1^2d\sigma_{S^{n-1}}^2,
$$
which can be considered as part of $S^{m+1}(\frac1{\sqrt{K_1}})\times_{(t_1+\frac{r_1}6)^{\gamma}}S^{n-1}$,
which obviously is of positive Ricci curvature.

Thus, the Ricci curvatures of $Q$ are positive.

\subsection{Curvatures of $P_i$}
For $P_i$, the metric is
$$
d\bar{s}^2=d\bar{t}^2+\bar{u}^2(\bar{t})d\theta_{S^m}^2+\bar{g}^2(\bar{t})d\sigma_{S^{n-1}}^2
$$
where $m\ge2, n\ge3$ and $n\ge m$ so that $2(n-1)\ge m\ge2$, and $0\le\bar t\le R_i$.
Similarly, let $T, \{\Theta_k\}, \{\Sigma_l\}$ be an orthonormal basis of the tangent space
corresponding to the directions $d\bar{t}, d\theta^2_{S^m}, d\sigma^2_{S^{n-1}}$ respectively.
Then, the sectional curvatures can be computed as follows.
\begin{align*}
K(T,\Theta_k,\Theta_k,T)&=-\frac{\bar{u}_{\bar{t}\bar{t}}}{\bar{u}},\\
K(\Theta_k,\Theta_p,\Theta_p,\Theta_k)&=\frac1{\bar{u}^2}-\frac{\bar{u}_{\bar{t}}^2}{\bar{u}^2},\\
K(T,\Sigma_l,\Sigma_l,T)&=-\frac{\bar{g}_{\bar{t}\bar{t}}}{\bar{g}},\\
K(\Sigma_l,\Sigma_q,\Sigma_q,\Sigma_l)&=\frac{1}{\bar{g}^2}-\frac{\bar{g}_{\bar{t}}^2}{\bar{g}^2},\\
K(\Theta_k,\Sigma_l,\Sigma_l,\Theta_k)&=-\frac{\bar{u}_{\bar{t}}}{\bar{u}}\frac{\bar{g}_{\bar{t}}}{\bar{g}}.
\end{align*}
And other terms of curvature tensors are zero.

The nonzero Ricci curvatures of $P_i$ are as follows.
\begin{align*}
Ric(T,T)&=-m\frac{\bar{u}_{\bar{t}\bar{t}}}{\bar{u}}-(n-1)\frac{\bar{g}_{\bar{t}\bar{t}}}{\bar{g}},\\
Ric(\Theta_k,\Theta_k)&=(m-1)(\frac1{\bar{u}^2}-\frac{\bar{u}_{\bar{t}}^2}{\bar{u}^2})
-\frac{\bar{u}_{\bar{t}\bar{t}}}{\bar{u}}
-(n-1)\frac{\bar{u}_{\bar{t}}}{\bar{u}}\frac{\bar{g}_{\bar{t}}}{\bar{g}},\\
Ric(\Sigma_l,\Sigma_l)&=(n-2)(\frac{1}{\bar{g}^2}-\frac{\bar{g}_{\bar{t}}^2}{\bar{g}^2})
-\frac{\bar{g}_{\bar{t}\bar{t}}}{\bar{g}}
-m\frac{\bar{u}_{\bar{t}}}{\bar{u}}\frac{\bar{g}_{\bar{t}}}{\bar{g}}.
\end{align*}

In the following, we will discuss the sectional and Ricci curvatures as $\bar t$
are in different intervals, i.e. $[0,d_i], [d_i,b_i]$, and $[d_i,R_i]$, respectively.

\subsubsection{In the interval $[0,d_i]$}

For $0\le\bar{t}\le d_i$,
\begin{align*}
K(T,\Theta_k,\Theta_k,T)&=-\frac{\bar{u}_{\bar{t}\bar{t}}}{\bar{u}}
\ge-(\frac1N+\frac1{N^2})l_i^2=O(\frac1{t_i^2}),\\
K(\Theta_k,\Theta_p,\Theta_p,\Theta_k)&=\frac1{\bar{u}^2}-\frac{\bar{u}_{\bar{t}}^2}{\bar{u}^2}
\ge\frac{1-\bar{u}_{\bar{t}^2(b_i)}}{\bar{u}^2}=(1-\frac1{N^2})\frac1{\bar{u}^2}
\ge(1-\frac1{N^2})l_i^2>0,\\
K(T,\Sigma_l,\Sigma_l,T)&=-\frac{\bar{g}_{\bar{t}\bar{t}}}{\bar{g}}=l_i^2>0,\\
K(\Sigma_l,\Sigma_q,\Sigma_q,\Sigma_l)&=\frac{1}{\bar{g}^2}-\frac{\bar{g}_{\bar{t}}^2}{\bar{g}^2}
=l_i^2>0,\\
K(\Theta_k,\Sigma_l,\Sigma_l,\Theta_k)&=-\frac{\bar{u}_{\bar{t}}}{\bar{u}}\frac{\bar{g}_{\bar{t}}}{\bar{g}}
=-\frac{l_i}{Nb_i}\frac{l_i\bar{t}}{\tan(l_i\bar{t})}\ge-\frac{l_i}{Nb_i}
\ge-\frac{l_i^2}{5n}=O(\frac1{t_i^2}).
\end{align*}
For the last estimate, note that it comes from $N\ge5n$ and $b_i\ge\frac1{l_i}$.
Thus, the condition of quadratically asymptotically nonnegative sectional curvatures holds in this interval.

For Ricci curvatures, we have
\begin{align*}
Ric(T,T)&=-m\frac{\bar{u}_{\bar{t}\bar{t}}}{\bar{u}}-(n-1)\frac{\bar{g}_{\bar{t}\bar{t}}}{\bar{g}}\\
        &\ge(n-1)l_i^2-m(\frac1N+\frac1{N^2})l_i^2\\
        &\ge\frac m2 l_i^2-m(\frac15+\frac1{25})l_i^2\\
        &>0,
\end{align*}
\begin{align*}
Ric(\Theta_k,\Theta_k)&=(m-1)(\frac1{\bar{u}^2}-\frac{\bar{u}_{\bar{t}}^2}{\bar{u}^2})
                        -\frac{\bar{u}_{\bar{t}\bar{t}}}{\bar{u}}
                        -(n-1)\frac{\bar{u}_{\bar{t}}}{\bar{u}}\frac{\bar{g}_{\bar{t}}}{\bar{g}}\\
                      &\ge(m-1)(1-\frac1{N^2})l_i^2-(\frac1N+\frac1{N^2})l_i^2-\frac{n-1}{5n}l_i^2\\
                      &\ge(1-\frac1{N^2}-\frac1N-\frac1{N^2}-\frac15)l_i^2\\
                      &>0,
\end{align*}
\begin{align*}
Ric(\Sigma_l,\Sigma_l)&=(n-2)(\frac{1}{\bar{g}^2}-\frac{\bar{g}_{\bar{t}}^2}{\bar{g}^2})
                       -\frac{\bar{g}_{\bar{t}\bar{t}}}{\bar{g}}
                       -m\frac{\bar{u}_{\bar{t}}}{\bar{u}}\frac{\bar{g}_{\bar{t}}}{\bar{g}}\\
                      &\ge(n-2)l_i^2+l_i^2-\frac m{5n}l_i^2\\
                      &\ge\big((n-1)-\frac15\big)l_i^2\\
                      &>0.
\end{align*}

\subsubsection{In the interval $[d_i,b_i]$}
For $d_i\le\bar{t}\le b_i$,
\begin{align*}
K(T,\Theta_k,\Theta_k,T)&=-\frac{\bar{u}_{\bar{t}\bar{t}}}{\bar{u}}
\ge-(\frac1N+\frac1{N^2})l_i^2=O(\frac1{t_i^2}),\\
K(\Theta_k,\Theta_p,\Theta_p,\Theta_k)&=\frac1{\bar{u}^2}-\frac{\bar{u}_{\bar{t}}^2}{\bar{u}^2}
\ge\frac{1-\bar{u}_{\bar{t}^2(b_i)}}{\bar{u}^2}=(1-\frac1{N^2})\frac1{\bar{u}^2}
\ge(1-\frac1{N^2})l_i^2>0,\\
K(T,\Sigma_l,\Sigma_l,T)&=-\frac{\bar{g}_{\bar{t}\bar{t}}}{\bar{g}}
=\frac{\gamma(1-\gamma)A}{\bar{t}^{2-\gamma}(A\bar{t}^\gamma+B)}>0,\\
K(\Sigma_l,\Sigma_q,\Sigma_q,\Sigma_l)&=\frac{1}{\bar{g}^2}-\frac{\bar{g}_{\bar{t}}^2}{\bar{g}^2}
=\frac{1-(\frac{\gamma A}{\bar{t}^{1-\gamma}})^2}{\bar{g}^2}
\ge\frac{1-(\frac{\gamma A}{d_i^{1-\gamma}})^2}{\bar{g}^2}=\frac{\sin^2(l_id_i)}{\bar{g}^2}>0,\\
K(\Theta_k,\Sigma_l,\Sigma_l,\Theta_k)&=-\frac{\bar{u}_{\bar{t}}}{\bar{u}}\frac{\bar{g}_{\bar{t}}}{\bar{g}}
=-\frac{l_i}{Nb_i}\frac{\gamma A\bar{t}^\gamma}{A\bar{t}^\gamma+B}
=-\frac{l_i}{Nb_i}\frac{\gamma A}{A+B\bar{t}^{-\gamma}}\\
&\ge-\frac{l_i}{Nb_i}\frac{\gamma A}{A+B d_i^{-\gamma}}
=-\frac{l_i}{Nb_i}\frac{l_id_i}{\tan(l_id_i)}\ge-\frac{l_i}{Nb_i}
\ge-\frac{l_i^2}{5n}.
\end{align*}
For the last estimate, note that it comes from $N\ge5n$ and $b_i\ge\frac1{l_i}$. Thus, the condition of
quadratically asymptotically nonnegative sectional curvatures also holds in this interval.

For Ricci curvatures,
\begin{align*}
Ric(T,T)&=-m\frac{\bar{u}_{\bar{t}\bar{t}}}{\bar{u}}-(n-1)\frac{\bar{g}_{\bar{t}\bar{t}}}{\bar{g}}\\
        &\ge-m(\frac1N+\frac1{N^2})l_i^2+\frac m2 \frac{\gamma(1-\gamma)A}{\bar{t}^{2-\gamma}(A\bar{t}^\gamma+B)}\\
        &\ge-m(\frac1N+\frac1{N^2})l_i^2+\frac m2 \frac{\gamma(1-\gamma)A}{b_i^{2-\gamma}(Ab_i^\gamma+B)}\\
        &\ge-m(\frac1N+\frac1{N^2})l_i^2+\frac m2 \frac{\gamma(1-\gamma)}{b_i^2}\\
        &\ge-m(\frac1N+\frac1{N^2})l_i^2+m\frac{\gamma(1-\gamma)}8l_i^2\\
        &>0,
\end{align*}
provided with $N\ge N(\gamma)$.

\begin{align*}
Ric(\Theta_k,\Theta_k)&=(m-1)(\frac1{\bar{u}^2}-\frac{\bar{u}_{\bar{t}}^2}{\bar{u}^2})
                        -\frac{\bar{u}_{\bar{t}\bar{t}}}{\bar{u}}
                        -(n-1)\frac{\bar{u}_{\bar{t}}}{\bar{u}}\frac{\bar{g}_{\bar{t}}}{\bar{g}}\\
                      &\ge(m-1)(1-\frac1{N^2})l_i^2-(\frac1N+\frac1{N^2})l_i^2-\frac{n-1}{5n}l_i^2\\
                      &\ge(1-\frac1{N^2}-\frac1N-\frac1{N^2}-\frac15)l_i^2\\
                      &>0,
\end{align*}
\begin{align*}
Ric(\Sigma_l,\Sigma_l)&=(n-2)(\frac{1}{\bar{g}^2}-\frac{\bar{g}_{\bar{t}}^2}{\bar{g}^2})
                       -\frac{\bar{g}_{\bar{t}\bar{t}}}{\bar{g}}
                       -m\frac{\bar{u}_{\bar{t}}}{\bar{u}}\frac{\bar{g}_{\bar{t}}}{\bar{g}}\\
                      &\ge(n-2)\frac{\sin^2(l_id_i)}{\bar{g}^2}+
                      \frac{\gamma(1-\gamma)A\bar{t}^\gamma}{\bar{t}^2(A\bar{t}^\gamma+B)}-
                      m\frac{l_i}{Nb_i}\frac{\gamma A\bar{t}^\gamma}{A\bar{t}^\gamma+B}\\
                      &>\frac{\gamma A\bar{t}^\gamma}{A\bar{t}^\gamma+B}
                      (\frac{1-\gamma}{\bar{t}^2}-m\frac{l_i}{Nb_i})\\
                      &\ge\frac{\gamma A\bar{t}^\gamma}{A\bar{t}^\gamma+B}
                      (\frac{1-\gamma}{b_i^2}-\frac{2m}{Nb_i^2})\\
                      &>0.
\end{align*}

\subsubsection{In the interval $[b_i,R_i]$}
For $b_i\le\bar{t}\le R_i$, by the previous construction, we have
\begin{align*}
0<\frac{\bar{u}_{\bar{t}}}{\bar{u}}&=\frac{1}{(1-C_3)\bar{t}+\frac{N}{l_i}-(1-C_3)b_i}\\
&\le\frac{1}{(1-C_3)\bar{t}+\frac{N}{l_i}-\frac{2(1-C_3)}{l_i}}\\
&\le\frac1{(1-C_3)\bar{t}}.
\end{align*}
Thus,
\begin{align*}
K(T,\Theta_k,\Theta_k,T)&=-\frac{\bar{u}_{\bar{t}\bar{t}}}{\bar{u}}
=-C_3\frac{\bar{u}_{\bar{t}}^2}{\bar{u}^2}\ge-\frac{C_3}{(1-C_3)^2\bar{t}^2}
\ge-\frac{C_3}{(1-C_3)^2b_i^2}=O(\frac1{t_i^2}),\\
K(\Theta_k,\Theta_p,\Theta_p,\Theta_k)&=\frac1{\bar{u}^2}-\frac{\bar{u}_{\bar{t}}^2}{\bar{u}^2}
\ge\frac{1-C_2^2}{\bar{u}^2}>0,\\
K(T,\Sigma_l,\Sigma_l,T)&=-\frac{\bar{g}_{\bar{t}\bar{t}}}{\bar{g}}
=\frac{\gamma(1-\gamma)A}{\bar{t}^{2-\gamma}(A\bar{t}^\gamma+B)}
\ge\frac{\gamma(1-\gamma)}{\bar{t}^2}\ge2\frac{C_3}{(1-C_3)^2\bar{t}^2}>0,
\end{align*}
provided with $C_3\le C_3(\gamma)$.
$$
K(\Sigma_l,\Sigma_q,\Sigma_q,\Sigma_l)=\frac{1}{\bar{g}^2}-\frac{\bar{g}_{\bar{t}}^2}{\bar{g}^2}
=\frac{1-(\frac{\gamma A}{\bar{t}^{1-\gamma}})^2}{\bar{g}^2}
\ge\frac{1-(\frac{\gamma A}{b_i^{1-\gamma}})^2}{\bar{g}^2}>0,
$$
provided with $t_i\ge t_i(c,r,\gamma,C_3,N)$. While
$$
K(\Theta_k,\Sigma_l,\Sigma_l,\Theta_k)=-\frac{\bar{u}_{\bar{t}}}{\bar{u}}\frac{\bar{g}_{\bar{t}}}{\bar{g}}
\ge-\frac{1-(1+C_3)C_2^2}{n-1}l_i^2=O(\frac1{t_i^2});
$$
Not that the above estimate will be verified after the calculation of $Ric(\Theta_k,\Theta_k)$ in the following.

Thus, the condition of quadratically asymptotically nonnegative sectional curvatures again holds in this interval.

For Ricci curvature,
\begin{align*}
Ric(T,T)&=-m\frac{\bar{u}_{\bar{t}\bar{t}}}{\bar{u}}-(n-1)\frac{\bar{g}_{\bar{t}\bar{t}}}{\bar{g}}\\
        &\ge-m\frac{C_3}{(1-C_3)^2\bar{t}^2}+2(n-1)\frac{C_3}{(1-C_3)^2\bar{t}^2}\\
        &>0,
\end{align*}
\begin{align*}
Ric(\Theta_k,\Theta_k)&=(m-1)(\frac1{\bar{u}^2}-\frac{\bar{u}_{\bar{t}}^2}{\bar{u}^2})
                        -\frac{\bar{u}_{\bar{t}\bar{t}}}{\bar{u}}
                        -(n-1)\frac{\bar{u}_{\bar{t}}}{\bar{u}}\frac{\bar{g}_{\bar{t}}}{\bar{g}}\\
                      &=(m-1)(\frac1{\bar{u}^2}-\frac{\bar{u}_{\bar{t}}^2}{\bar{u}^2})
                        -C_3\frac{\bar{u}_{\bar{t}}^2}{\bar{u}^2}
                        -(n-1)\frac{\bar{u}_{\bar{t}}}{\bar{u}}\frac{\bar{g}_{\bar{t}}}{\bar{g}}\\
                      &\ge\frac{1-(1+C_3)C_2^2}{C_2^2}
                      \frac{(\frac{t_i}{C_1C_2})^{{2C_3}/(1-C_3)}}
                      {\left((1-C_3)\bar{t}+\frac{N}{l_i}-(1-C_3)b_i\right)^{{2}/(1-C_3)}}\\
                      &\quad-\frac{n-1}{(1-C_3)\bar{t}+\frac{N}{l_i}-(1-C_3)b_i}
                      \frac{\gamma A\bar{t}^{\gamma-1}}{A\bar{t}^\gamma+B}\\
                      &\ge\frac{1-(1+C_3)C_2^2}{C_2^2}
                      \frac{(\frac{t_i}{C_1C_2})^{{2C_3}/(1-C_3)}}
                      {\left((1-C_3)\bar{t}+\frac{N}{l_i}-(1-C_3)b_i\right)^{2/(1-C_3)}}\\
                      &\quad-\frac{n-1}{(1-C_3)\bar{t}+\frac{N}{l_i}-(1-C_3)b_i}
                      \frac{\gamma A\bar{t}^{\gamma-1}}{\frac12A\bar{t}^\gamma}\\
                      &=\frac1{\bar{t}\big((1-C_3)\bar{t}+\frac{N}{l_i}-(1-C_3)b_i\big)^{2/(1-C_3)}}
                      \Big(\frac{1-(1+C_3)C_2^2}{C_2^2}(\frac{t_i}{C_1C_2})^{\frac{2C_3}{1-C_3}}\bar{t}\\
                      &\quad-2(n-1)\gamma\big((1-C_3)\bar{t}+\frac{N}{l_i}-(1-C_3)b_i\big)^{\frac{1+C_3}{1-C_3}}\Big).
\end{align*}
Let $H(\bar{t})\triangleq\frac{1-(1+C_3)C_2^2}{C_2^2}(\frac{t_i}{C_1C_2})^{\frac{2C_3}{1-C_3}}\bar{t}-
2(n-1)\gamma\left((1-C_3)\bar{t}+\frac{N}{l_i}-(1-C_3)b_i\right)^{\frac{1+C_3}{1-C_3}}$, $b_i\le\bar{t}\le R_i$.
Note that
\begin{align*}
H(b_i)&=\frac{1-(1+C_3)C_2^2}{C_2^2}(\frac{t_i}{C_1C_2})^{\frac{2C_3}{1-C_3}}b_i-
      2(n-1)\gamma(\frac N{l_i})^{\frac{1+C_3}{1-C_3}}\\
      &\ge\frac{1-(1+C_3)C_2^2}{C_2^2}(\frac{t_i}{C_1C_2})^{\frac{2C_3}{1-C_3}}\frac1{l_i}-
      2(n-1)\gamma(\frac N{l_i})^{\frac{1+C_3}{1-C_3}}\\
      &=\frac{C(c,r,C_3)}{N^{\frac1{C_3}}}t_i^{\frac{1+C_3}{1-C_3}}-
      \frac{2(n-1)\gamma C'(c,r,C_3)}{N^{1+\frac1{C_3}}}t_i^{\frac{1+C_3}{1-C_3}}\\
      &>0,
\end{align*}
provided with $N>N(c,r,\gamma,C_3)$.

On the other hand,
\begin{align*}
H_{\bar{t}}(\bar{t})&=\frac{1-(1+C_3)C_2^2}{C_2^2}(\frac{t_i}{C_1C_2})^{\frac{2C_3}{1-C_3}}-
                    2(n-1)\gamma(1+C_3)\big((1-C_3)\bar{t}\\
                    &\quad+\frac{N}{l_i}-(1-C_3)b_i\big)^{\frac{2C_3}{1-C_3}}\\
                    &\ge\frac{1-(1+C_3)C_2^2}{C_2^2}(\frac{t_i}{C_1C_2})^{\frac{2C_3}{1-C_3}}-
                    2(n-1)\gamma(1+C_3)(\frac{t_i}{C_1C_2})^{\frac{2C_3}{1-C_3}}\\
                    &=\left(\frac1{C_2^2}-(1+C_3)-2(n-1)\gamma(1+C_3)\right)(\frac{t_i}{C_1C_2})^{\frac{2C_3}{1-C_3}}\\
                    &\ge0,
\end{align*}
provided with $C_3\le\frac12(\frac1{C_2^2}-1)$ and $\gamma\le\frac1{8(n-1)}(\frac1{C_2^2}-1)$.

Thus, for $b_i\le\bar{t}\le R_i$,
$$
H(\bar{t})\ge H(b_i)>0.
$$
So, $Ric(\Theta_k,\Theta_k)>0$. In addition, we also get
$$
-\frac{\bar{u}_{\bar{t}}}{\bar{u}}\frac{\bar{g}_{\bar{t}}}{\bar{g}}
>-\frac{1-(1+C_3)C_2^2}{n-1}\bar{u}^2\ge-\frac{1-(1+C_3)C_2^2}{n-1}l_i^2,
$$
which gives the required estimate of $K(\Theta_k, \Sigma_l, \Sigma_l, \Theta_k)$ before. Thus,
\begin{align*}
Ric(\Sigma_l,\Sigma_l)&=(n-2)(\frac{1}{\bar{g}^2}-\frac{\bar{g}_{\bar{t}}^2}{\bar{g}^2})
                       -\frac{\bar{g}_{\bar{t}\bar{t}}}{\bar{g}}
                       -m\frac{\bar{u}_{\bar{t}}}{\bar{u}}\frac{\bar{g}_{\bar{t}}}{\bar{g}}\\
                      &\ge(n-2)\frac{1-(\gamma Ab_i^{\gamma-1})^2}{\bar{g}^2}+
                      \frac{\gamma(1-\gamma)A\bar{t}^\gamma}{\bar{t}^2(A\bar{t}^\gamma+B)}-
                      \frac m{n-1}\left(1-(1+C_3)C_2^2\right)l_i^2\\
                      &>\frac1{2\bar{g}^2}-2\left(1-(1+C_3)C_2^2\right)l_i^2\\
                      &\ge\frac1{2(1+\frac6r)^{2\gamma}t_i^{2\gamma}}-2\left(1-(1+C_3)C_2^2\right)l_i^2\\
                      &>0,
\end{align*}
provided with $t_i$ sufficiently large.

\section{Gluing along the boundary and smoothing near the boundary}
In this section, we will first check that $Q\setminus(B_{\frac45r_i}(o_i)\times_{g_i}S^{n-1})$ and $P_i$ have isometric boundaries so that we can glue them together along the boundaries to get a $C^0$ metric. Next, by means of the normal curvature properties of the boundaries, we are able to construct a $C^1$ metric to replace the above $C^0$ metric near the glued boundaries and to verify that the quadratically asymptotic non-negativeness of curvature and positiveness of Ricci curvature are still preserved.

\subsection{Gluing along the boundary}
First, we consider the boundary of $Q\setminus(B_{\frac45r_i}(o_i)\times_{g_i}S^{n-1})$, that is $\partial B_{\frac45r_i}(o_i)\times_{g_i}S^{n-1}$. By the previous construction, $[t_i+\frac{r_i}6,t_i+\frac{11r_i}6]\times_{u(t)}\times S^m$ can be considered as part of $S^{m+1}(\frac1{\sqrt{K_i}})$; so, the removed geodesic ball $B_{\frac45r_i}(o_i)$ is a geodesic ball of radius $\frac45r_i$ in $S^{m+1}(\frac1{\sqrt{K_i}})$. On the other hand, for $t_i+\frac{r_i}6<t<t_i+\frac{11r_i}6$, $g(t)\equiv g_i=(t_i+\frac {r_i}6)^\gamma=(1+\frac r6)^\gamma t_i^\gamma$. So, after a suitable rotation, the restricted metric on the boundary $\partial B_{\frac45r_i}(o_i)\times_{g_i}S^{n-1}$ is
$$
(\frac{\sin\frac45\sqrt{K}r}{\sqrt{K}}t_i)^2d\theta_{S^m}^2+\left((1+\frac r6)^\gamma t_i^\gamma\right)^2d\sigma_{S^{n-1}}^2.
$$

Next, we consider the boundary of $P_i$, which is exactly the hypersurface $\bar{t}=R_i$ with the restricted metric
$$
(\frac{\sin{\frac45}\sqrt{K}r}{\sqrt{K}}t_i)^2d\theta_{S^m}^2+\big((1+\frac r6)^\gamma t_i^\gamma\big)^2d\sigma_{S^{n-1}}^2.
$$
Thus, both boundaries are isometric.

\subsection{Smoothing near the boundary}
As seen above, a $C^0$ metric $h_0$ is already constructed on the manifold $M^{m+n}$, i.e
$$
\left(M^{m+n},h_0\right)=\left(Q\setminus\coprod\limits_{i=1}^{+\infty}(B_{\frac45r_i}(o_i)\times_{g_i}S^{n-1}),ds^2\right)
\cup_{\textrm{Id}}\coprod\limits_{i=1}^{+\infty}\left(P_i,d\bar{s}_i^2\right);
$$
$h_0$ is smooth on $M^{m+n}$ except the glued boundary part $\partial B_{\frac45r_i}(o_i)\times_{g_i}S^{n-1}$ on which $h_0$ is $C^0$ and has positive Ricci curvature and quadratically asymptotically nonnegative sectional curvature on the smooth part. Furthermore, the restricted metric of $h_0$ on the $2\varepsilon_i$-neighbourhood ($\varepsilon_i\ll1\ll r_i$) of the boundary (i.e. $h_0\mid_{[-2\varepsilon_i,2\varepsilon_i]\times S^m\times S^{n-1}}$) can be rewritten as
\begin{align*}
\left\{
      \begin{array}{ll}
      dt^2+\left(\frac{\sin\sqrt{K_i}(t+\frac45r_i)}{\sqrt{K_i}}\right)^2d\theta_{S^m}^2+\left((t_i+\frac {r_i}6)^\gamma\right)^2d\sigma_{S^{n-1}}^2, & 0\le t \le2\varepsilon_i, \vspace{+0.2cm}\\
      dt^2+\bar{u}^2(t+R_i)d\theta_{S^m}^2+\bar{g}^2(t+R_i)d\sigma_{S^{n-1}}^2, & -2\varepsilon_i\le t <0,
      \end{array}
\right.
\end{align*}
which, as mentioned above, is $C^0$ at $t=0$ since
$$
\bar{u}(R_i)=\frac{\sin\frac45\sqrt{K}r}{\sqrt{K_i}},\quad \bar{g}(R_i)=g_i=(t_i+\frac {r_i}6)^\gamma.
$$

In the following, we will construct a $C^1$ metric $h_1$ on $[-2\varepsilon_i,2\varepsilon_i]\times S^m\times S^{n-1}$, which is the same as $h_0$ on $[-2\varepsilon_i,-\varepsilon_i]\times S^m\times S^{n-1}$ and $[\varepsilon_i,2\varepsilon_i]\times S^m\times S^{n-1}$, and verify that it also has positive Ricci curvature and quadratically asymptotically nonnegative curvature on this part (for similar construction also see \cite{BWW}, but our case will be more complicated due to the condition of quadratically asymptotic non-negative curvature).

As before, we will introduce some constants for convenience,
$$
a=\frac{\sin\sqrt{K_i}(\frac45r_i+\varepsilon_i)}{\sqrt{K_i}},\quad b=\cos\sqrt{K_i}(\frac45r_i+\varepsilon_i),
$$
$$
c=\bar{u}(R_i-\varepsilon_i),\quad d=\bar{u}_{\bar{t}}(R_i-\varepsilon_i),
$$
$$
e=\bar{g}(R_i-\varepsilon_i),\quad f=\bar{g}_{\bar{t}}(R_i-\varepsilon_i).
$$

The $C^1$ metric then is
\begin{align*}
h_1=\left\{
      \begin{array}{ll}
      dt^2+\left(\frac{\sin\sqrt{K_i}(t+\frac45r_i)}{\sqrt{K_i}}\right)^2d\theta_{S^m}^2+\left((t_i+\frac {r_i}6)^\gamma\right)^2d\sigma_{S^{n-1}}^2, & \varepsilon_i\le t \le2\varepsilon_i, \vspace{+0.2cm}\\
      dt^2+U^2(t)d\theta_{S^m}^2+G^2(t)d\sigma_{S^{n-1}}^2, & -\varepsilon_i<t<\varepsilon_i, \vspace{+0.2cm}\\
      dt^2+\bar{u}^2(t+R_i)d\theta_{S^m}^2+\bar{g}^2(t+R_i)d\sigma_{S^{n-1}}^2, & -2\varepsilon_i\le t\le-\varepsilon_i,
      \end{array}
\right.
\end{align*}
where
\begin{align*}
U^2(t)&=\frac{(c^2-a^2)+2(ab+cd)\varepsilon_i}{4\varepsilon_i^3}t^3+\frac{ab-cd}{2\varepsilon_i}t^2 -\frac{3(c^2-a^2)+2(ab+cd)\varepsilon_i}{4\varepsilon_i}t \\
      &\quad+\frac{(a^2+c^2)+(cd-ab)\varepsilon_i}2,\\
G^2(t)&=\frac{(e^2-g_i^2)+2ef\varepsilon_i}{4\varepsilon_i^3}t^3-\frac{ef}{2\varepsilon_i}t^2 -\frac{3(e^2-g_i^2)+2ef\varepsilon_i}{4\varepsilon_i}t+\frac{(e^2+g_i^2)+ef\varepsilon_i}2.
\end{align*}
It is clear that $h_1$ is the same as $h_0$ on $[-2\varepsilon_i,-\varepsilon_i]\times S^m\times S^{n-1}$ and $[\varepsilon_i,2\varepsilon_i]\times S^m\times S^{n-1}$ respectively.

By a direct computation, we have
\begin{align*}
U_tU&=\frac{3(c^2-a^2)+6(ab+cd)\varepsilon_i}{8\varepsilon_i^3}t^2+\frac{ab-cd}{2\varepsilon_i}t -\frac{3(c^2-a^2)+2(ab+cd)\varepsilon_i}{8\varepsilon_i}, \\
U_t^2+U_{tt}U&=\frac{3(c^2-a^2)+6(ab+cd)\varepsilon_i}{4\varepsilon_i^3}t+\frac{ab-cd}{2\varepsilon_i},
\end{align*}
and
\begin{align*}
G_tG&=\frac{3(e^2-g_i^2)+6ef\varepsilon_i}{8\varepsilon_i^3}t^2-\frac{ef}{2\varepsilon_i}t -\frac{3(e^2-g_i^2)+2ef\varepsilon_i}{8\varepsilon_i}, \\
G_t^2+G_{tt}G&=\frac{3(e^2-g_i^2)+6ef\varepsilon_i}{4\varepsilon_i^3}t-\frac{ef}{2\varepsilon_i}.
\end{align*}

Puting $t=\pm\varepsilon_i$ into the above equations gives
$$
U(\varepsilon_i)=a=\frac{\sin\sqrt{K_i}(\frac45r_i+\varepsilon_i)}{\sqrt{K_i}},\quad U_t(\varepsilon_i)=b=\cos\sqrt{K_i}(\frac45r_i+\varepsilon_i),
$$
$$
U(-\varepsilon_i)=c=\bar{u}(R_i-\varepsilon_i),\quad U_t(-\varepsilon_i)=d=\bar{u}_{\bar{t}}(R_i-\varepsilon_i),
$$
$$
G(\varepsilon_i)=g_i,\quad G_t(\varepsilon_i)=0,
$$
$$
G(-\varepsilon_i)=e=\bar{g}(R_i-\varepsilon_i),\quad G_t(-\varepsilon_i)=f=\bar{g}_{\bar{t}}(R_i-\varepsilon_i).
$$
Thus, the metric $h_1$ on $[-2\varepsilon_i,2\varepsilon_i]\times S^m\times S^{n-1}$ is smooth for $t\neq\pm\varepsilon_i$ and $C^1$ at $t=\pm\varepsilon_i$. Using this $C^1$ metric $h_1\mid_{[-2\varepsilon_i,2\varepsilon_i]\times S^m\times S^{n-1}}$ to replace the $C^0$ metric $h_0\mid_{[-2\varepsilon_i,2\varepsilon_i]\times S^m\times S^{n-1}}$, we obtain a new $C^1$ metric on the whole manifold $M^{m+n}$, still denoted by $h_1$, which actually is smooth except for $t=\pm\varepsilon_i$.

Next, we will verify that the new metric $h_1$ still has positive Ricci curvature and quadratically asymptotically nonnegative curvature. For this purpose, we only need to consider the curvature terms of
$$
dt^2+U^2(t)d\theta_{S^m}^2+G^2(t)d\sigma_{S^{n-1}}^2,\quad -\varepsilon_i<t<\varepsilon_i.
$$
Similar to what we have done, we have
\begin{align*}
K(T,\Theta_k,\Theta_k,T)&=-\frac{U_{tt}}{U},\\
K(\Theta_k,\Theta_p,\Theta_p,\Theta_k)&=\frac1{U^2}-\frac{U_t^2}{U^2},\\
K(T,\Sigma_l,\Sigma_l,T)&=-\frac{G_{tt}}{G},\\
K(\Sigma_l,\Sigma_q,\Sigma_q,\Sigma_l)&=\frac{1}{G^2}-\frac{G_t^2}{G^2},\\
K(\Theta_k,\Sigma_l,\Sigma_l,\Theta_k)&=-\frac{U_t}{U}\frac{G_t}{G}.
\end{align*}
And other terms of curvature tensors are zero.

Note that when $\varepsilon_i\ll1\ll r_i$,
$$
c^2-a^2=-2\frac{\sin\frac45\sqrt{K}r}{\sqrt{K_i}}\left(\cos\frac35\sqrt{K}r+\cos\frac45\sqrt{K}r\right)\varepsilon_i+ o(\varepsilon_i),
$$
$$
ab+cd=\frac{\sin\frac45\sqrt{K}r}{\sqrt{K_i}}\left(\cos\frac35\sqrt{K}r+\cos\frac45\sqrt{K}r\right)+o(1),
$$
$$
a^2+c^2=2\left(\frac{\sin\frac45\sqrt{K}r}{\sqrt{K_i}}\right)^2+o(1),
$$
$$
ab-cd=-\frac{\sin\frac45\sqrt{K}r}{\sqrt{K_i}}\left(\cos\frac35\sqrt{K}r-\cos\frac45\sqrt{K}r\right)+o(1),
$$
$$
e^2-g_i^2=-2g_i\bar{g}_{\bar{t}}(R_i)\varepsilon_i+ o(\varepsilon_i),
$$
$$
ef=g_i\bar{g}_{\bar{t}}(R_i)+o(1),
$$
$$
e^2+g_i^2=2g_i^2+o(1).
$$
Then
\begin{align*}
U^2(t)&=\left(\frac{\sin\frac45\sqrt{K}r}{\sqrt{K_i}}\right)^2+o(1), \\
U_tU&=\frac12\frac{\sin\frac45\sqrt{K}r}{\sqrt{K_i}}\left(\cos\frac35\sqrt{K}r+\cos\frac45\sqrt{K}r\right) \\ &\quad-\frac12\frac{\sin\frac45\sqrt{K}r}{\sqrt{K_i}}\left(\cos\frac35\sqrt{K}r-\cos\frac45\sqrt{K}r\right) \frac{t}{\varepsilon_i}+o(1), \\
U_t^2+U_{tt}U&=-\frac12\frac{\sin\frac45\sqrt{K}r}{\sqrt{K_i}}\left(\cos\frac35\sqrt{K}r-\cos\frac45\sqrt{K}r\right) \cdot\frac1{\varepsilon_i}+O(1),
\end{align*}
and
\begin{align*}
G^2(t)&=g_i^2+o(1), \\
G_tG&=\frac12g_i\bar{g}_{\bar{t}}(R_i)-\frac12g_i\bar{g}_{\bar{t}}(R_i)\frac{t}{\varepsilon_i}+o(1), \\
G_t^2+G_{tt}G&=-\frac12g_i\bar{g}_{\bar{t}}(R_i)\cdot\frac{1}{\varepsilon_i}+O(1),
\end{align*}
which together imply that
$$
-\frac{U_{tt}}{U}=\frac{\frac12\frac{\sin\frac45\sqrt{K}r}{\sqrt{K_i}}\left(\cos\frac35\sqrt{K}r -\cos\frac45\sqrt{K}r\right)\cdot\frac1{\varepsilon_i}+O(1) +U_t^2}{\left(\frac{\sin\frac45\sqrt{K}r}{\sqrt{K_i}}\right)^2+o(1)}>0,
$$
$$
-\frac{G_{tt}}{G}=\frac{\frac12g_i\bar{g}_{\bar{t}}(R_i)\cdot\frac{1}{\varepsilon_i}+O(1)+G_t^2}{g_i^2+o(1)}>0,
$$
\begin{align*}
\left|\frac{U_t}{U}\right|&=\frac{\frac12\frac{\sin\frac45\sqrt{K}r}{\sqrt{K_i}}\left[\left(\cos\frac35\sqrt{K}r+ \cos\frac45\sqrt{K}r\right) -\left(\cos\frac35\sqrt{K}r- \cos\frac45\sqrt{K}r\right)\frac{t}{\varepsilon_i}\right] +o(1)}{\left(\frac{\sin\frac45\sqrt{K}r}{\sqrt{K_i}}\right)^2+o(1)} \\
                          &\le\frac{\cos\frac35\sqrt{K}r+o(1)}{\frac{\sin\frac45\sqrt{K}r}{\sqrt{K_i}}+o(1)} =\frac{\bar{u}_{\bar{t}}(R_i)+o(1)}{\bar{u}(R_i)+o(1)},
\end{align*}
and
$$
\left|\frac{G_t}{G}\right|=\left|\frac{\frac12g_i\bar{g}_{\bar{t}}(R_i)\left(1-\frac{t}{\varepsilon_i}\right) +o(1)}{g_i^2+o(1)}\right|\le\frac{\bar{g}_{\bar{t}}(R_i)+o(1)}{\bar{g}(R_i)+o(1)}.
$$
By choosing $\varepsilon_i$ sufficiently small, $-\frac{U_{tt}}{U}$ and $-\frac{G_{tt}}{G}$ can have arbitrarily large lower bound, meanwhile the upper bound of $\left|\frac{U_t}{U}\right|$ and $\left|\frac{G_t}{G}\right|$ can be arbitrarily close to $\frac{\bar{u}_{\bar{t}}(R_i)}{\bar{u}(R_i)}$ and $\frac{\bar{g}_{\bar{t}}(R_i)}{\bar{g}(R_i)}$ respectively. Based on these fact, both the quadratically asymptotic non-negativeness of curvature and the positiveness of Ricci curvature will be remained after replacing the $C^0$ metric by the $C^1$ metric.

\vskip .2cm
\noindent {\bf Remark:} As mentioned in the introduction, we should remark that in the $C^1$ smoothing from the $C^0$ metric near the glued boundaries, the positiveness of Ricci curvature (equivalently, positive lower bounds of curvature terms $-\frac{U_{tt}}{U}$ and $-\frac{G_{tt}}{G}$) is controlled by the cofficients "$\cos\frac35\sqrt{K}r -\cos\frac45\sqrt{K}r$" and "$\bar{g}_{\bar{t}}(R_i)$" of the terms with $\frac{1}{\varepsilon_i}$. In fact, the positiveness of these terms is just the condition in Perelman's gluing criterion, i.e. "the normal curvatures of $\partial M_1$ are bigger than the negative of the normal curvatures of $\partial M_2$", that is why we can construct an explicit $C^1$ smoothing and preserve the corresponding curvature conditions in our case.

\vskip.5cm
\noindent Mathematics \& Science College, Shanghai Normal University, Shanghai, 200234\\
{\it E-mail}: jianghuihong@shnu.edu.cn
\vskip .3cm
\noindent Department of Mathematics, Shanghai Jiao Tong University, Shanghai, 200240\\
{\it E-mail}: {yangyihu@sjtu.edu.cn}

\end{document}